\documentclass[11pt]{amsart}

\usepackage{amssymb}
\usepackage{amsmath}
\usepackage{amsthm}
\usepackage{color} 
\usepackage{verbatim}
\usepackage{bbm}
\usepackage{dsfont}

\newtheorem{theorem}{Theorem}[section]

\newtheorem{remark}{Remark}[section]

\numberwithin{equation}{section}

\newcommand{\Rd}{\mathbb{R}^d}

\newcommand{\vx}{v_{(\xi)}}
\newcommand{\vxx}{v_{(\xi)(\xi)}}

\newcommand{\px}{\psi_{(\xi)}}

\newcommand{\pxsop}{\frac{\psi_{(\xi)}^2}{\psi}}

\begin{document}

\title[Regularity of degenerate Monge-Amp\`ere equations]{Interior regularity of fully nonlinear degenerate elliptic equations, II: real and complex Monge-Amp\`ere equations}
\author{Wei Zhou}
\address{School of Mathematics, University of Minnesota}
\email{zhoux123@math.umn.edu}


\begin{abstract}
\noindent
We first consider the Dirichlet problem for the degenerate real Monge-Amp\`ere equation:
\begin{equation}\tag{Real M-A}\label{rm-a}
\left\{\begin{array}{rcll}
\det(u_{x^ix^j})&=&f(x)  &\text{in }D\qquad\qquad\ \\ 
u&=&g  &\text{on }\partial D.
\end{array}
\right.
\end{equation}
We prove that if $D\subset \Rd\ (d\ge2)$ is a bounded strictly convex domain with $C^{3}$ boundary, $g\in C^{1,1}(\partial D)$, $0\le f^{1/d}\in C^{0,1}(\bar D)$ and there exists a constant $K$ such that $f^{1/d}+K|x|^2$ is convex in $\bar D$, then there exists a convex function $u\in C^{1,1}_{loc}(D)\cap C^{0,1}(\bar D)$ uniquely solving (\ref{rm-a}).

Then we consider the Dirichlet problem for the degeneratete complex Monge-Amp\`ere equation:
\begin{equation}\label{cm-a}\tag{Complex M-A}
\left\{\begin{array}{rcll}
\det(u_{z^j\bar z^k})&=&f(z)  &\text{in }D\qquad\qquad\qquad\\ 
u&=&g  &\text{on }\partial D.
\end{array}
\right.
\end{equation}
We prove that if $D\subset \mathbb C^d$ ($d\ge2$) is a bounded strictly pseudoconvex domain with $C^{3}$ boundary, $g\in C^{1,1}(\partial D)$ and $0\le f^{1/d}\in C^{1,1}(\bar D)$, then there exists a plurisubharmonic function $u\in C^{1,1}_{loc}(D)\cap C^{0,1}(\bar D)$ uniquely solving (\ref{cm-a}). 


We also estimate the derivatives up to second order in both problems. Our technique is probabilistic by following Krylov's approach.

\end{abstract}

\maketitle

\section{Introduction}


The first problem we study is the degenerate real Monge-Amp\`ere equation with Dirichlet boundary condition:
\begin{equation}\label{RMA}\tag{\bf{Real M-A}}
\left\{\begin{array}{rcll}
\det(u_{x^ix^j})&=&f  &\text{in }D\qquad\qquad\qquad\\ 
u&=&g  &\text{on }\partial D,
\end{array}
\right.
\end{equation}
where the domain $D\subset \Rd$ is bounded, strictly convex and sufficiently smooth,
$(u_{x^ix^j})_{d\times d}$ is the Hessian matrix of $u$ and the function $f=f(x)$ is nonnegative. We seek a convex function uniquely solving (\ref{RMA}) and investigate its regularity. This problem is very important in many fields and might be very challenging depending on the regularity and geometry of $D$, the regularity of $f$ and $g$, as well as the positivity property of $f$. It has been studied extensively by many people. Particularly, ($D$ is a bounded and strictly convex domain, unless specificed.)
\begin{itemize}
\item In \cite{MR0437805}, Cheng and Yau proved that if $D$ is $C^{2}$-smooth, $g\in C^2(\bar D)$, $f\in C_{loc}^\infty(D)$ and $0<f\le B \operatorname{dist}(x,\partial D)^{\beta-d-1}$, for some $B>0$, $\beta>0$, then $u\in C_{loc}^{\infty}(D)$; if $\beta=d+1$, then $u\in C^{0,1}(\bar D)$. They also obtained that if $D$ is convex (not necessarily strictly convex), $g=0$, $f\in C_{loc}^{\infty}(D)$ and $f>0$, then $u\in C^0(\bar D)\cap C^{\infty}_{loc}(D)$.
\item In \cite{MR688919}, Krylov showed that if $D$ is of $C^\infty$, $g\in C^{\infty}(\partial D)$, $f\in C^{\infty}(\bar D)$ and $f>0$, then $u \in C^{\infty}(\bar D)$.
\item In \cite{MR739925}, Caffarelli, Nirenberg and Spruck obtained the same result, and they also showed in \cite{MR864651} that if $D$ is $C^{3,1}$-smooth, $g\in C^{3,1}(\partial D)$ and $f=0$, then $u\in C^{1,1}(\bar D)$.
\item In \cite{MR766792}, Trudinger and Urbas proved that if $D$ is $C^{1,1}$-smooth, $f=0$ and $g\in C^{1,1}(\bar D)$, then $u\in C^{0,1}(\bar D)\cap C^{1,1}_{loc}(D)$. 
\item In \cite{MR992979}, Krylov showed that if $D$ is $C^{3,1}$-smooth, $g\in C^{3,1}(\partial D)$, $f^{1/d}\in C^{1,1}(\bar D)$ and $f\ge0$, then $u\in C^{1,1}(\bar D)$.
\item In \cite{MR1687172}, Guan, Trudinger and Wang obtained that if $D$ is $C^{3,1}$-smooth, $g\in C^{3,1}(\partial D)$, $f^{1/(d-1)}\in C^{1,1}(\bar D)$ and $f\ge0$, then $u\in C^{1,1}(\bar D)$, which is optimal in the sense of the regularity assumption on $f$, due to an example by Wang in \cite{MR1223269}. 
\end{itemize}

For the degenerate Monge-Amp\`ere equation, we first note that $C^{1,1}$-regularity is the best that we can expect, even if the boundary data $g$ is analytic on $\partial D$. This can be seen by an example given in \cite{MR1687172}, by considering the unit ball in $\mathbb R^2$ as the domain $D$ and
$$u(x_1,x_2)=\big[\max\{(x_1^2-1/2)^+,(x_2^2-1/2)^+\}\big]^2.$$
We also note that the assumption that $g\in C^{3,1}(\partial D)$ is necessary for obtaining the global $C^{1,1}$-regularity of $u$. See, for example, Example 1 in \cite{MR864651}. While the sufficiency of $g\in C^{3,1}(\partial D)$ to obtain the global $C^{1,1}$-regularity is established by the aforementioned papers \cite{MR864651, MR992979, MR1687172} under various settings. Therefore, an interesting problem is investigating the interior $C^{1,1}$-regularity of the solution to (\ref{RMA}), when g is only assumed to be in the class of $C^{1,1}(\partial D)$, whose necessity is obvious. The sufficiency for the homogeneous case was obtained by Trudinger and Urbas in the above mentioned paper \cite{MR766792}. Trudinger and Urbas's proof of this result relied on the fact that in this homogeneous case the Monge-Amp\`ere equation reduces to $\lambda_{\operatorname{min}}(u_{xx})=0$, which implies that the epigraph of the solution is the convex envelope of the epigraph of the boundary data $\varphi$. Our result generalizes theirs in the sense of considering $f\ge0$ in general, rather than $f\equiv0$.

By thinking of the Monge-Amp\`ere equation as a special Bellman equation with constant coefficients, we obtain regularity and solvability  results on the degenerate real Monge-Amp\`ere equation with Dirichlet boundary condition. We assume $D$ be a bounded, $C^{3}$-smooth and strictly convex domain and $f$ be a function from $D$ to $[0,\infty)$. Based on Krylov's viewpoint, we introduce the probabilistic solution to (\ref{RMA}):
$$v(x)=\inf_{\alpha\in\mathfrak A}E\bigg[g(x^{\alpha,x}_{\tau^{\alpha,x}})-\int_0^{\tau^{\alpha,x}}d\sqrt[d]{\det\big((1/2)\alpha_t\alpha_t^*\big)f(x_t^{\alpha,x})}dt\bigg],$$
with
$$x_t^{\alpha,x}=x+\int_0^t\alpha_sdw_s,$$
where $w_t$ is a Wiener process and $\mathfrak A$ is the set of progressively-measurable processes $\alpha_t$ with values in $\mathbb R^{d\times d}$ satisfying $\operatorname{tr}(\alpha_t\alpha_t^*)=2$ for all $t\ge0$.  Our main result for (\ref{RMA}) is as follows:
\begin{itemize}
\item If $f\in C^{0,1}(\bar D)$, $g\in C^{1,1}(\partial D)$, and there exists a positive constant $K$ such that $f^{1/d}+K|x|^2$ is convex in $\bar D$, then $v\in C^{1,1}_{loc}(D)\cap C^{0,1}(\bar D)$, and for a.e. $x\in D$, we have the second derivative estimate
$$0\le v_{(\xi)(\xi)}\le N\bigg(|\xi|^2+\pxsop\bigg),\ \forall\xi\in \Rd,$$
where $\psi$ is the defining function for $D$.
Meanwhile, $v$ is the unique convex solution to (\ref{RMA}) in the space of $C^{1,1}_{loc}(D)\cap C^{0,1}(\bar D)$.
\end{itemize}

The interior $C^{1,1}$-regularity is optimal under our regularity assumption on the boundary data $g$. The second derivative estimate coincides with the fact that the $C^{1,1}$-norm of $v$ shouldn't blow up along the tangent directions on the boundary, and says that the $C^{1,1}$-norm of $v$ doesn't blow up faster than $1/\operatorname{dist}(\cdot,\partial D)^2$ in any direction near the boundary. 


The second problem we study is the degenerate complex Monge-Amp\`ere equation with Dirichlet boundary condition:
\begin{equation}\label{CMA}\tag{\bf{Complex M-A}}
\left\{\begin{array}{rcll}
\det(u_{z^j\bar z^k})&=&f  &\text{in }D\qquad\qquad\qquad\qquad\\ 
u&=&g  &\text{on }\partial D,
\end{array}
\right.
\end{equation}
where the domain $D \subset \mathbb C^d$ is bounded, strictly pseudoconvex and sufficiently smooth, $(u_{z^j \bar z^k})_{d\times d}$ is the Hessian matrix of $u$ and the function $f=f(z)$ is nonnegative. We seek a plurisubharmonic function uniquely solving (\ref{CMA}) and investigate its regularity.

Compared with its real counterpart, the regularity theory on the degenerate complex Monge-Amp\`ere equation is much less developed. Many regularity results for the real Monge-Amp\`ere equation can not be extended to the complex case, because of various reasons. Important breakthroughs have been made by the following works: ($D$ is a bounded and strictly pseudoconvex domain.)

\begin{itemize}
\item In \cite{MR0445006}, Bedford and Taylor proved that if $D$ is the unit ball, $g\in C^2(\partial D)$, $f^{1/d}\in C^2(\bar D)$ and $f\ge0$, then $u\in C^{1,1}_{loc}(D)\cap C(\bar D)$. They also showed that when $D$ is smooth but not necessarily a ball, if $f^{1/d}\in C^{0,1}(\bar D)$ and $g\in C^{2}(\bar D)$, then $u\in C^{0,1}(\bar D)$.

\item In \cite{MR780073}, Caffarelli, Kohn, Nirenberg and Spruck showed that if $D$ is of $C^{\infty}$, $g\in C^{\infty}(\partial D)$, $f\in C^{\infty}(\bar D)$ and $f>0$, then $u\in C^{\infty}(\bar D)$. Regularity results under various conditions on $f=f(z,u,u_z)$ were also established there.

\item In \cite{MR992979}, Krylov obtained that if $D$ is $C^{3,1}$-smooth, $g\in C^{3,1}(\partial D)$, $f^{1/d}\in C^{1,1}(\bar D)$ and $f\ge0$, then $u\in C^{1,1}(\bar D)$. This seems to be the only known $C^{1,1}$-regularity result up to boundary for the degenerate case, even when the domain $D$ is the unit ball.
\end{itemize}

Bedford and Taylor's approach in \cite{MR0445006} made use of the transitivity of the automorphism group of the unit ball in $\mathbb C^d$, which is not applicable for general strictly pseudoconvex domains. Our result generalizes Bedford and Taylor's interior $C^{1,1}$-regularity result by allowing the domain $D$ be any bounded, $C^{3}$-smooth and strictly pseudoconvex domain. To be precise, we let $D\subset \mathbb C^d$ be a bounded, $C^{3}$-smooth and strictly pseudoconvex domain and $f$ be a function from $D$ to $[0,\infty)$. Define the probabilistic solution to (\ref{CMA}) by 
$$v(z)=\inf_{\alpha\in\mathfrak A}E\bigg[g(z^{\alpha,z}_{\tau^{\alpha,z}})-\int_0^{\tau^{\alpha,z}}d\sqrt[d]{\det\big(\alpha_t\bar\alpha_t^*\big)f(z_t^{\alpha,z})}dt\bigg],$$
with
$$z_t^{\alpha,x}=z+\int_0^t\alpha_sdW_s,$$
where $W_t$ is a normalized complex Wiener process
and $\mathfrak A$ is the set of progressively measurable processes $\alpha_t$ with values in $\mathbb C^{d\times d}$ satisfying $\operatorname{tr}(\alpha_t\bar\alpha_t^*)=1$ for all $t\ge0$.  Our main result for (\ref{CMA}) is the following:
\begin{itemize}
\item If $f^{1/d}\in C^{1,1}(\bar D)$ and $g\in C^{1,1}(\partial D)$, then $u\in C^{1,1}_{loc}(D)\cap C^{0,1}(\bar D)$, and for a.e. $z\in D$, we have the second derivative estimate
$$0\le v_{(\xi)(\xi)}\le N\bigg(|\xi|^2+\pxsop\bigg),\ \forall\xi\in \mathbb C^d,$$
where $\psi$ is the defining function for $D$. Furthermore, $v$ is the unique plurisubharmonic solution to (\ref{CMA}) in the space of $C^{1,1}_{loc}(D)\cap C^{0,1}(\bar D)$.
\end{itemize}

Again, the interior $C^{1,1}$-regularity is optimal under our regularity assumption on the boundary data $g$. The second derivative estimate says that the $C^{1,1}$-norm of $v$ doesn't blow up faster than $1/\operatorname{dist}(\cdot,\partial D)^2$ in any non-tangent direction near the boundary. 



This paper is outlined as follows. In Section 2, we deal with the degenerate real Monge-Amp\`re equation by directly applying our results in \cite{InteriorRegularityI}. In Section 3, we extend our results in Section 2 to the degenerate complex Monge-Amp\`ere equation. 

We end this section by introducing the notation. Throughout the article, the summation convention for repeated indices is assumed.
Given any sufficiently smooth function $u(x)$ from $\Rd$ to $\mathbb R$, for $y,z\in\Rd$, let
\begin{align*}
u_{(y)}=& u_{x^i}y^i,\qquad  u_{(y)(z)}= u_{x^i x^j}y^i z^j,\qquad u_{(y)}^2=(u_{(y)})^2.
\end{align*}
We denote the gradient vector of $u$ by $u_x$ and the Hessian matrix of $u$ by $u_{xx}$.
For any matrix $\sigma$, its tranpose is denoted by $\sigma^*$.

The notation of complex analysis will be introduced at the beginning of Section 3.

\section{Interior $C^{1,1}$ regularity of the degenerate real Monge-Amp\`ere equation}
In this section, we consider the Dirichlet problem for the degenerate real Monge-Amp\`ere equation in a strictly convex domain.

\subsection{Statement of the theorem}
Let $D$ be a bounded domain in $\Rd$ described by a $C^{3}$ function $\psi$ which is non-singular on $\partial D$, i.e.
$$D:=\{x\in \Rd:\psi(x)>0\},\qquad|\psi_x|\ge 1\mbox{ on }\partial D.$$
We also assume that $\psi$ is strictly concave in $\bar D$, i.e.
\begin{equation}\label{scr}
\forall a\in \bar{\mathcal S}^+_d: \operatorname{tr}(a)=1,\qquad\operatorname{tr}(a\psi_{xx})<0\mbox{ in }\bar D,
\end{equation}
where $\bar{\mathcal S}^+_d$ denotes the set of all non-negative symmetric $d\times d$ matrices.

Let $w_t$ be a Wiener process of dimension $d$, and $\mathfrak A$ be the set of progressively-measurable processes $\alpha_t$ with values in $\mathbb R^{d\times d}$ satisfying $\operatorname{tr}(\alpha_t\alpha_t^*)=2, \forall t\ge0$. Introduce a family of controlled diffusion processes
$$x_t^{\alpha,x}=x+\int_0^t\alpha_sdw_s,  \qquad \forall \alpha_t\in\mathfrak A.$$
Denote $\tau^{\alpha,x}$ the first exit time of $x_t^{\alpha,x}$ from $D$.

Let $f$ and $g$ be bounded measurable functions on $\bar D$ with values in $[0,\infty)$ and $\mathbb R$, respectively.

\begin{theorem}\label{real}
Let
\begin{equation}\label{prma}
v(x)=\sup_{\alpha\in\mathfrak A}E\bigg[g(x^{\alpha,x}_{\tau^{\alpha,x}})+\int_0^{\tau^{\alpha,x}}\sqrt[d]{\det(a^\alpha_t)}f(x^{\alpha,x}_t)dt\bigg], 
\end{equation}
with 
$$a^\alpha_t=\frac{1}{2}\alpha_t\alpha_t^*.$$

If $f,g\in C^{0,1}(\bar D)$, then $v\in C^{0,1}_{loc}(D)\cap C(\bar D)$, and for a.e. $x\in D$,
\begin{equation}\label{firstma}
|\vx|\le N\bigg(|\xi|+\frac{|\px|}{\psi^{1/2}}\bigg),\ \forall\xi\in \Rd,
\end{equation}
where the constant $N=N(|f|_{0,1,D},|g|_{0,1,D},|\psi|_{3,D},d)$.

If $f\in C^{0,1}(\bar D)$, $g\in C^{1,1}(\bar{D})$, and there exists a constant $K$ such that $f+K|x|^2$ is convex in $\bar D$, then $v\in C^{1,1}_{loc}(D)\cap C^{0,1}(\bar D)$, and for a.e. $x\in D$, 
\begin{equation}\label{secondma}
-N\bigg(|\xi|^2+\pxsop\bigg)\le \vxx\le 0,\ \forall\xi\in\Rd,
\end{equation}
where the constant $N=N(|f|_{0,1,D},|g|_{1,1,D}, |\psi|_{3,D},K, d)$. Meanwhile, $u=-v$ is the unique convex solution in $C^{1,1}_{loc}(D)\cap C^{0,1}(\bar D)$ of the Dirichlet problem for the degenerate Monge-Amp\`ere equation:
\begin{equation}\label{realma}
\left\{\begin{array}{rcll}
\det(u_{xx})&=&d^{-d}f^d  &\text{a.e. in }D\\ 
u&=&-g  &\text{on }\partial D,
\end{array}
\right.
\end{equation}
satisfying the second derivative estimate: for a.e. $x\in D$,
\begin{equation*}
0\le u_{(\xi)(\xi)}\le N\bigg(|\xi|^2+\pxsop\bigg),\ \forall\xi\in \Rd.
\end{equation*}
\end{theorem}

\begin{remark}
The assumption that $f+K|x|^2$ is convex in $\bar D$ for some constant $K$ is weaker than the assumption that $f\in C^{1,1}(\bar D)$.
\end{remark}

\begin{remark}
The first derivative estimate (\ref{firstma}) is not new. It has been obtained in \cite{MR1211724}, in which an example (see Remark 5.1 in \cite{MR1211724}) was also given showing that (\ref{firstma}) is sharp. The interior $C^{1,1}$-regularity result and the second derivative estimate (\ref{secondma}) are new. However, the author doesn't know whether the estimate is sharp.
\end{remark}

\subsection{Proof of Theorem \ref{real}}
We first show that Theorems 2.1, 2.2 and 2.3 in \cite{InteriorRegularityI} are applicable for $v(x)$ given in (\ref{prma}) by verifying Assumptions 2.1 and 2.2 in \cite{InteriorRegularityI} and the weak non-degeneracy condition (See Remark 2.1 in \cite{InteriorRegularityI}). Then we show that the associated dynamic programming equation of $v(x)$ is equivalent to the Monge-Amp\`ere equation in (\ref{realma}).

\begin{proof}
We apply Theorems 2.1, 2.2 and 2.3 in \cite{InteriorRegularityI} with
$$A=\{\alpha\in\mathbb R^{d\times d}:\operatorname{tr}(\alpha\alpha^*)=2\},$$
$$\sigma^\alpha=\alpha,\qquad a^\alpha=\frac{1}{2}\alpha\alpha^*,\qquad b^\alpha=0,\qquad c^\alpha=0,\qquad f^\alpha=\sqrt[d]{\det(a^\alpha)}f.$$
In this situation,
$$\mathbb A=\{a^\alpha:\alpha\in A\}\subset\{a\in\bar{\mathcal S}^+_d: \operatorname{tr}(a)=1\}.$$
It follows that (\ref{scr}) reads
$$\sup_{\alpha\in A}\operatorname{tr}(a^\alpha\psi)<0.$$
Due to the compactness of $\mathbb A$, by replacing $\psi$ with $N\psi$, we may assume
$$\sup_{\alpha\in A}\operatorname{tr}(a^\alpha\psi)\le-1,$$
which is exactly Assumption 2.1 in \cite{InteriorRegularityI}.

Next, we notice that for any orthogonal matrix $O$ of size $d\times d$,
$$\operatorname{tr}(Oa^\alpha O^*)=\operatorname{tr}(O^*Oa^\alpha)=\operatorname{tr}(a^\alpha)=1,$$
which implies that
$$O\mathbb A O^*\subset\mathbb A=O^*(O\mathbb A O^*)O\subset O\mathbb A O^*.$$
Therefore Assumption 2.2 in \cite{InteriorRegularityI} also holds. (See also Remark 2.1 in \cite{InteriorRegularityI}.)

To verify to weak non-degeneracy we consider $a^{\alpha_0}=(1/d)I_{d\times d}\in\mathbb A$, then by (2.5) and (2.6) in \cite{InteriorRegularityI} we have
$$\mu=\inf_{|\zeta|=1}\sup_{\alpha\in A}(a^\alpha)_{ij}\zeta^i\zeta^j\ge\inf_{|\zeta|=1}(a^{\alpha_0})_{ij}\zeta^i\zeta^j=1/d>0.$$
Therefore Theorems 2.1, 2.2 and 2.3 in \cite{InteriorRegularityI} are true for $v(x)$ defined by (\ref{prma}). 

It remains to prove the second inequality in (\ref{secondma}) and the equivalence between the associated dynamic programming equation 
\begin{equation}\label{bellmanma}
\sup_{\alpha\in A}\Big[(a^\alpha)_{ij}v_{x^ix^j}+\sqrt[d]{\det(a^\alpha)}f\Big]=0
\end{equation}
and the Monge-Amp\`ere equation in (\ref{realma}). 
They are well-known facts, which were actually Lemma 2 in Section 3.2 of \cite{MR901759}. For the sake of completeness and the convenience of the reader we give the following argument.

First, we rewrite (\ref{bellmanma}) as
\begin{equation}\label{bellmanmaa}
\sup_{a\in\mathbb A}\Big[\operatorname{tr}(av_{xx})+\sqrt[d]{\det(a)}f\Big]=0.
\end{equation}

In particular, in $D$,
\begin{equation}\label{asdf}
\operatorname{tr}(av_{xx})+\sqrt[d]{\det(a)}f\le0
\end{equation}
for each $a\in\mathbb A$. For any fixed $\zeta\in\Rd$ with $|\zeta|=1$ and $a=\zeta\zeta^*$, from (\ref{asdf}) we get
$$(v_{xx}\zeta,\zeta)\le -\sqrt[d]{\det(a)}f\le0,$$
which proves the second inequality in (\ref{secondma}).

Next, take $\delta>0$ and set
$$a=(\delta I-v_{xx})^{-1}c_\delta,\qquad c_\delta^{-1}=\operatorname{tr}\big[(\delta I-v_{xx})^{-1}\big].$$
Then $\operatorname{tr}(a)=1$ and (\ref{asdf}) yields
$$\operatorname{tr}\big[(\delta I-v_{xx})^{-1}v_{xx}\big]+\big[\det(\delta I-v_{xx})\big]^{-1/d}f\le0.$$
It follows that
\begin{align*}
f\le&\big[\det(\delta I-v_{xx})\big]^{1/d}\operatorname{tr}\big[(\delta I-v_{xx})^{-1}(-v_{xx})\big]\\
=&\big[\det(\delta I-v_{xx})\big]^{1/d}\big\{\operatorname{tr}(I)-\delta\operatorname{tr}\big[(\delta I- v_{xx})^{-1}\big]\big\}\\
\le&\big[\det(\delta I-v_{xx})\big]^{1/d}d.
\end{align*}
Therefore we have
$$d^{-d}f^d\le\det(\delta I-v_{xx}).$$
By letting $\delta\downarrow0$ we obtain
$$d^{-d}f^d\le\det(-v_{xx})=\det(u_{xx}).$$
In fact, here an equality holds instead of the inequality. To prove this suppose that at some point $x_0\in D$, we have
$$d^{-d}f^d<\det(u_{xx}).$$
Then, in particular, $\det(u_{xx})>0$ and $u_{xx}$ is non-degenerate at $x_0$. Take a matrix $a_0\in \mathbb A$ which attains the supremum in (\ref{bellmanmaa}) at $x_0$. Since
$$-\operatorname{tr}(a_0v_{xx})=\operatorname{tr}(a_0u_{xx})\ge \lambda_{\min}\operatorname{tr}(a_0)>0,$$
where $\lambda_{\min}$ is the smallest eigenvalue of the strictly positive matrix $u_{xx}$, we see that
$$\det(a_0)>0.$$
Now by the fact that the geometric mean is not bigger than the arithmetic mean, 
\begin{align*}
f\big[\det(a_0)\big]^{1/d}=&-\operatorname{tr}(a_0v_{xx})=\operatorname{tr}(a_0u_{xx})=\operatorname{tr}\big(\sqrt{a_0}u_{xx}\sqrt{a_0}\big)\\
\ge& d\big[\det\big(\sqrt{a_0}u_{xx}\sqrt{a_0}\big)\big]^{1/d}=d\big[\det(u_{xx})\big]^{1/d}\big[\det(a_0)\big]^{1/d}.
\end{align*}
It follows that at $x_0$,
$$d^{-d}f^d\ge\det(u_{xx}),$$
which gives the desired contradiction. Hence the equivalence of (\ref{bellmanma}) and the Monge-Amp\`ere equation in (\ref{realma}) is proved.

\end{proof}
\section{Interior $C^{1,1}$ regularity of the degenerate complex Monge-Amp\`ere equation}

In this section, we consider the Dirichlet problem for degenerate complex Monge-Amp\`ere equation in a strictly pseudoconvex domain.
\subsection{Statement of the theorem}
We use the following standard notation: $\mathbb C$ denotes the set of all complex numbers; and
\begin{align*}
z=(z^1,\cdots,z^d)=&(x^1+ix^{d+1},\cdots,x^k+ix^{d+k},\cdots, x^d+ix^{d+d})\\
=&(x^1,\cdots,x^d)+i(x^{d+1},\cdots,x^{d+d})=:\operatorname{Re}z+i\operatorname{Im}z
\end{align*}
is an element of $\mathbb C^d$.  We also use the following notation of partial differential operators:
$$u_{z^k}=\frac{1}{2}(u_{x^k}-iu_{y^k}),\qquad u_{\bar z^k}=\frac{1}{2}(u_{x^k}+iu_{y^k})$$
$$u_{z^k\bar z^j}=(u_{z^k})_{\bar z^j}\qquad u_{z\bar z}=\big(u_{z^k\bar z^j}\big)_{1\le j,k\le d}$$
Moreover, for any $\xi,\eta\in \mathbb C^d$, we define
$$u_{(\xi)}=u_{z^k}\xi^k+u_{\bar z^k}\bar\xi^k,\qquad u_{(\xi)(\eta)}=\big(u_{(\xi)}\big)_{(\eta)}$$
Since any function $u$ from $\mathbb C^d$ to $\mathbb R$ can be viewed as a function from $\mathbb R^{2d}$ to $\mathbb R$, by abuse of notation we write $u(z)=u(x)$ with $x=(\operatorname{Re}z,\operatorname{Im}z)$. As a result, we see that
$$u_{(\xi)}(z)=u_{(\operatorname{Re\xi},\operatorname{Im}\xi)}(x),\qquad u_{(\xi)(\eta)}(z)=u_{(\operatorname{Re\xi},\operatorname{Im}\xi)(\operatorname{Re\eta},\operatorname{Im}\eta)}(x)$$

Let $D$ be a bounded domain in $\mathbb C^d$ described by a $C^{3}$ function $\psi$ which is non-singular on $\partial D$, i.e.
$$D:=\{z\in \mathbb C^d:\psi(z)>0\},\qquad|\psi_z|\ge 1\mbox{ on }\partial D.$$
We also assume that $\psi$ is strictly plurisuperharmonic in $\bar D$, i.e.
\begin{equation}\label{sc}
\forall a\in \bar{\mathcal H}^+_d: \operatorname{tr}(a)=1,\qquad\operatorname{tr}(a\psi_{z\bar z})<0\mbox{ in }\bar D,
\end{equation}
where $\bar{\mathcal H}^+_d$ denotes the set of all non-negative Hermitian $d\times d$ matrices.

Let $W_t$ be a normalized complex Wiener process of dimension $d$, i.e. a $d$-dimensional stochastic process $W_t=(W_t^1,..., W_t^d)$ with values in $\mathbb C^d$ given by 
$$W_t^j=\frac{1}{\sqrt 2}\Big(w_t^{j,1}+iw_t^{j,2}\Big), \ t\ge0, \ 1\le j\le d,$$
where the processes $(w_t^{j,1},w_t^{j,2})_{1\le j\le d}$ are independent real Wiener processes.

Let $\mathfrak A$ be the set of progressively-measurable processes $\alpha_t$ with values in $\mathbb C^{d\times d}$ satisfying $\operatorname{tr}(\alpha_t\bar\alpha_t^*)=1, \forall t\ge0$, and
$$A=\{\alpha\in \mathbb C^{d\times d}: \operatorname{tr}(\alpha\bar\alpha^*)=1\}.$$
If we denote $\alpha\bar\alpha^*$ by $a^\alpha$, then
$$\mathbb A:=\{a^\alpha:\alpha\in A\}=\{a\in\bar{\mathcal H}_d^+: \operatorname{tr}a=1\}.$$
Introduce a family of controlled complex diffusion processes
$$z_t^{\alpha,x}=z+\int_0^t\alpha_sdW_s,  \qquad \forall \alpha_t\in\mathfrak A.$$
Denote $\tau^{\alpha,z}$ the first exit time of $z_t^{\alpha,z}$ from $D$.

Let $f$ and $g$ be bounded measurable functions on $\bar D$ with values in $[0,\infty)$ and $\mathbb R$ respectively.

\begin{theorem}\label{complex}
Let
\begin{equation}\label{prmac}
v(z)=\sup_{\alpha\in\mathfrak A}E\bigg[g(z^{\alpha,z}_{\tau^{\alpha,z}})+\int_0^{\tau^{\alpha,z}}\sqrt[d]{\det(a^{\alpha}_t)}f(z^{\alpha,z}_t)dt\bigg], 
\end{equation}
with
$$ a^\alpha_t=\alpha_t\bar\alpha_t^*.$$

If $f,g\in C^{0,1}(\bar D)$, then $v\in C^{0,1}_{loc}(D)\cap C(\bar D)$, and for a.e. $z\in D$,
\begin{equation}\label{firstmac}
|\vx|\le N\bigg(|\xi|+\frac{|\px|}{\psi^{1/2}}\bigg), \ \forall \xi\in\mathbb C^d,
\end{equation}
where the constant $N=N(|f|_{0,1,D},|g|_{0,1,D},|\psi|_{3,D},d)$.

If $f, g\in C^{1,1}(\bar{D})$, then $v\in C^{1,1}_{loc}(D)\cap C^{0,1}(\bar D)$, and for a.e. $z\in D$, 
\begin{equation}\label{secondmac}
-N\bigg(|\xi|^2+\pxsop\bigg)\le \vxx\le 0,\ \forall \xi\in\mathbb C^d,
\end{equation}
where the constant $N=N(|f|_{1,1,D},|g|_{1,1,D},|\psi|_{3,D}, d)$. Meanwhile, $u=-v$ is the unique plurisubharmonic solution in $C^{1,1}_{loc}(D)\cap C^{0,1}(\bar D)$ of the Dirichlet problem for the degenerate complex Monge-Amp\`ere equation:
\begin{equation}\label{complexma}
\left\{\begin{array}{rcll}
\det(u_{z\bar z})&=&d^{-d}f^d  &\text{a.e. in }D\\ 
u&=&-g  &\text{on }\partial D,
\end{array}
\right.
\end{equation}
satisfying the second derivative estimate: for a.e. $x\in D$,
\begin{equation*}
0\le u_{(\xi)(\xi)}\le N\bigg(|\xi|^2+\pxsop\bigg),\ \forall\xi\in\mathbb C^d.
\end{equation*}
\end{theorem}

\begin{remark}
When studying the $C^{1,1}$-regularity of $v$, the regularity assumption on $f$ we need  is actually weaker than $C^{1,1}(\bar D)$. Similarly to the real case, we only need that the generalized second derivatives of $f$ are bounded from below. Therefore, the assumption that $f\in C^{1,1}(\bar D)$ can be replace with the assumptions that $f+K|x|^2$ is convex for some constant $K$, where we treat $f=f(x)$ as a function of $2d$ real variables with $x=(\operatorname{Re}z,\operatorname{Im}z)\in \mathbb R^{2d}$. When $f$ is sufficiently smooth, the condition that $f+K|x|^2$ is convex is equivalent to
$$\operatorname{Re}\Big(f_{z^jz^k}(z)\xi^j\xi^k\Big)+f_{z^j\bar z^k}(z)\xi^j\bar\xi^k\ge-K|\xi|^2,\ \forall \xi\in \mathbb C^{d}.$$

\end{remark}

\subsection{Proof of Theorem \ref{complex}}
We prove Theorem \ref{complex} by making use of Theorems 2.1, 2.2 and 2.3 in \cite{InteriorRegularityI}.

\begin{proof}
We first define the following homomorphisms:
$$\Phi: \mathbb C^d\rightarrow\mathbb R^{2d}; z\mapsto 
\begin{pmatrix}
\operatorname{Re}z\\
\operatorname{Im}z
\end{pmatrix}
$$
and
$$
\Phi: \mathbb C^{d\times d}\rightarrow\mathbb R^{2d\times 2d}; \alpha\mapsto 
\begin{pmatrix}
\operatorname{Re}\alpha&\operatorname{Im}\alpha\\
-\operatorname{Im}\alpha&\operatorname{Re}\alpha
\end{pmatrix}.
$$
To rewrite the value function in (\ref{prmac}) as a function on $\mathbb R^{2d\times 2d}$, we notice that
$$\Phi( z_t^{\alpha,z})=\Phi z+\int_0^t\frac{1}{\sqrt 2}(\Phi\alpha_s)dw_s$$
where $w_t$ is a Wiener process of dimension $2d$.

For any Hermitian matrix $a$, there is a unitary matrix $U$ such that $Ua\bar U^*$ is a real diagonal matrix $M$. 
\begin{align*}
\det(\Phi a)=&\det\big[\Phi(\bar U^*MU)\big]=\det\big[\Phi(\bar U^*)\Phi M\Phi U)\big]\\
=&\det\big[(\Phi U)^*\Phi M\Phi U)\big]=\det(\Phi M)=\big[\det(M)\big]^2.
\end{align*}
Therefore
$$\sqrt[d]{\det(a^\alpha_t)}=\sqrt[2d]{\det\Phi (\alpha_t\bar\alpha_t^*)}=\sqrt[2d]{\det(\Phi \alpha_t\Phi\bar\alpha_t^*)}=\sqrt[2d]{\det\big[(\Phi \alpha_t)(\Phi\alpha_t)^*\big]}.$$
If we use the notation:
$$\beta=\Phi\alpha, \qquad\mathfrak B=\Phi\mathfrak A=\{\Phi\alpha:\alpha\in\mathfrak A\},$$
then we can rewrite (\ref{prmac}) as
\begin{equation}
v(x)=\sup_{\beta\in\mathfrak B}E_x^{\beta}\bigg[g(x_{\tau})+\int_0^{\tau}\sqrt[2d]{\det(\beta_t\beta_t^*)}f(x_t)dt\bigg],
\end{equation}
where
$$x_t^{\beta,x}=x+\int_0^t\frac{1}{\sqrt 2}\beta_sdw_s.$$
By noticing that
$$\operatorname{tr}\big[(\Phi a)\psi_{xx}\big]=4\operatorname{tr}(a\psi_{z\bar z}),$$
Assumption 2.1 in \cite{InteriorRegularityI} is satisfied by replacing $\psi$ with $N\psi$ for sufficiently constant $N$.
However, if we have a try on applying Theorems 2.1-2.3 in \cite{InteriorRegularityI}  directly, we should fail at Assupmtion 2.2 in \cite{InteriorRegularityI}. Because Assupmtion 2.2 in \cite{InteriorRegularityI} doesn't hold for
$$\mathbb B=\Phi \mathbb A=\{\Phi (a^\alpha):a^\alpha\in\mathbb A\}.$$
More precisely, since
$$\mathbb B=\bigg\{
\begin{pmatrix}
S&T\\
-T&S
\end{pmatrix}
: S\in\bar{\mathcal S}^+_d, T \mbox{ is skew symmetric}, \operatorname{tr}(S)=1\bigg\}
,$$
the relation
$$O\mathbb B O^*=\mathbb B$$
doesn't hold for all orthogonal matrix $O$ of size $2d\times 2d$.

Fortunately, for any unitary matrix of size $d\times d$, we have
$$U\mathbb A\bar U^*=\mathbb A.$$
which can play a role of Assumption 2.2 in \cite{InteriorRegularityI}. Indeed, we observe that
$$\mathbb B=\Phi\mathbb A=\Phi(U\mathbb A\bar U^*)=\Phi U\Phi\mathbb A\Phi(\bar U^*)=(\Phi U)\mathbb B(\Phi U)^*.$$
Moreover, we note that if $Q\in \mathbb C^{d\times d}$ is skew Hermitian, then $e^Q$ is unitary, and therefore
$e^{\Phi Q}=\Phi(e^Q)$ satisfies
$$(e^{\Phi Q})\mathbb B (e^{\Phi Q})^*=\mathbb B$$
Therefore, in order to apply Theorems 2.1-2.3 in \cite{InteriorRegularityI}, it suffice to find a suitable matrix function $P=P(x,\xi)$ from $D\times\mathbb R^{2d}$ to $\mathbb R^{2d\times2d}$, such that $P$ can be expressed as $\Phi Q$, where $Q$ is skew Hermitian, and $P(x,\xi)$ satisfies all properties it has in the proof of Lemma 7.1 in \cite{InteriorRegularityI}. 

To construct $P$ let us start from looking at the equation (7.1) in \cite{InteriorRegularityI}, the most crucial property it satisfies in the proof of Lemma 7.1 in \cite{InteriorRegularityI}.
We define
$$\chi: D_\delta^\lambda\times \mathbb C^d\rightarrow \mathbb C; (z,\xi)\mapsto -\frac{\psi_{\bar z^k}(\psi_{z^k})_{(\xi)}}{|\psi_{\bar z}|^2}$$
and
$$R: D_\delta^\lambda\times \mathbb C^d\rightarrow \mathbb C^{d\times d}; (z,\xi)\mapsto \bigg(\frac{(\psi_{z^k})_{(\xi)}\psi_{\bar z^j} -(\psi_{\bar z^j})_{(\xi)}\psi_{z^k}}{|\psi_{\bar z}|^2}\bigg)_{1\le j,k\le d},$$
which are analogous to $\rho$ and $P$ in Lemma 7.1 in \cite{InteriorRegularityI}.

We claim that
\begin{equation}\label{nm}
(\psi_{\bar z})_{(\xi)}+R\psi_{\bar z}+\chi\psi_{\bar z}=0.
\end{equation}
Indeed, we have
\begin{align*}
&\Big[(\psi_{\bar z})_{(\xi)}+R\psi_{\bar z}+\chi\psi_{\bar z}\Big]^j\\
=&(\psi_{\bar z^j})_{(\xi)}+\frac{(\psi_{ z^k})_{(\xi)}\psi_{\bar z^j}-(\psi_{\bar z^j})_{(\xi)}\psi_{z^k}}{\psi_z\psi_{\bar z}}\psi_{\bar z^k}-\frac{\psi_{\bar z^k}(\psi_{ z^k})_{(\xi)}}{\psi_z\psi_{\bar z}}\psi_{\bar z^j}\\
=&(\psi_{\bar z^j})_{(\xi)}\Big[1-\frac{\psi_{z^k}\psi_{\bar z^k}}{\psi_z\psi_{\bar z}}\Big]=0.
\end{align*}
Next, we notice that $\chi$ is not real in general, so we decompose it as $\chi=\rho+i\varkappa$, where $\rho$ and $\varkappa$ are real valued. If we denote $R+i\varkappa I$ as $Q$, the equation (\ref{nm}) can be rewritten as
\begin{equation}\label{important}
(\psi_{\bar z})_{(\xi)}+Q\psi_{\bar z}+\rho\psi_{\bar z}=0.
\end{equation}
We emphasis that $\rho$ is real and $Q$ is skew Hermitian.
From (\ref{important}) and the fact that $\psi_x=2\Phi(\psi_{\bar z})$, we obtain
\begin{align*}
0=&\Big(\alpha\alpha^*\big((\psi_{\bar z})_{(\xi)}+Q\psi_{\bar z}+\rho\psi_{\bar z}\big),(\psi_{\bar z})_{(\xi)}+Q\psi_{\bar z}+\rho\psi_{\bar z}\Big)\\
=&\bigg(\Phi\Big(\alpha\alpha^*\big((\psi_{\bar z})_{(\xi)}+Q\psi_{\bar z}+\rho\psi_{\bar z}\big)\Big),\Phi\Big((\psi_{\bar z})_{(\xi)}+Q\psi_{\bar z}+\rho\psi_{\bar z}\Big)\bigg)\\
=&\frac{1}{4}\Big(\beta\beta^*\big((\psi_x)_{(\xi)}+(\Phi Q)\psi_x+\rho\psi_x \big),(\psi_x)_{(\xi)}+(\Phi Q)\psi_x+\rho\psi_x\Big)
\end{align*}
Therefore if we let $P=\Phi Q$,
$$\psi_{(\xi)(\beta^k)}+\rho\psi_{(\beta^k)}+\psi_{(P\beta^k)}=0.$$
It is also not hard to see that $P$ and $\rho$ we define here satisfy all the other property in Lemma 7.1 in \cite{InteriorRegularityI}. As a result, we can apply Theorem 2.1-2.3 in \cite{InteriorRegularityI} to obtain all regularity results of $v$ defined by (\ref{prmac}) stated in Theorem \ref{complex}.

It remains to verify that the associated dynamic programing equation is equivalent to the complex Monge-Amp\`ere equation in (\ref{complexma}). To write down the real Bellman equation we note that its diffusion term is $(1/\sqrt 2)\beta$. Hence the associated dynamic programing equation is the real Bellman equation
\begin{equation}
\sup_{\beta\in B}\Big[(1/4)\operatorname{tr}(\beta\beta^*v_{xx})+\sqrt[2d]{\det(\beta\beta^*)}f\Big]=0,
\end{equation}
which is equivalent to 
\begin{equation}
\sup_{a\in A}\Big\{\operatorname{tr}\big[(\Phi a)v_{xx})\big]+4\sqrt[2d]{\det(\Phi a)}f\Big\}=0.
\end{equation}
To write down the corresponding complex Bellman equation, it suffices to notice that
$$\operatorname{tr}\big[(\Phi a)v_{xx}\big]=4\operatorname{tr}(av_{z\bar z}).$$
Therefore the complex Bellman equation is
\begin{equation}\label{bellmanmab}
\sup_{a\in A}\Big[\operatorname{tr}(av_{z\bar z})+\sqrt[d]{\det( a)}f\Big]=0,
\end{equation}
which has the same form of (\ref{bellmanma}). Therefore, the equivalence between (\ref{bellmanmab}) and the complex Monge-Amp\`ere equation in (\ref{complexma}) can be verified by repeating the argument right after (\ref{bellmanmaa}). The proof is complete.
\end{proof}

\section{Acknowledgements} The author wishes to express sincere gratitude towards his PhD advisor, Professor Nicolai V. Krylov, for illuminating suggestions and the financial support during the preparation of this paper. The author is also very grateful to Professor Hongjie Dong for inspiring discussions on fully nonlinear PDE theory.


\providecommand{\bysame}{\leavevmode\hbox to3em{\hrulefill}\thinspace}
\providecommand{\MR}{\relax\ifhmode\unskip\space\fi MR }
\providecommand{\MRhref}[2]{%
  \href{http://www.ams.org/mathscinet-getitem?mr=#1}{#2}
}
\providecommand{\href}[2]{#2}

\end{document}